\documentclass[12pt]{article}
\usepackage{amscd}
\usepackage{amsfonts}
\usepackage{amsmath,amsthm,hyperref}
\usepackage[all]{xy}
\usepackage{amsmath,amsthm,hyperref}
\usepackage{amsmath,amssymb,amsthm,latexsym}
\usepackage{amscd}
\newtheorem{theorem}{Theorem}[section]
\newtheorem{proposition}[theorem]{Proposition}
\newtheorem{definition}[theorem]{Definition}
\newtheorem{corollary}[theorem]{Corollary}
\newtheorem{lemma}[theorem]{Lemma}

\theoremstyle{remark}
\newtheorem{remark}[theorem]{Remark}

\begin{document}
\title{On Willmore surfaces in $S^n$ of flat normal bundle}
\author{Peng Wang\footnote{ Supported by Program for Young Excellent Talents in Tongji University and the Tianyuan Foundation of China, grant 10926112.} }
\date{}
\maketitle

\begin{center}
{\bf Abstract}\\
\end{center}
We discuss several kinds of Willmore surfaces of flat normal bundle in this paper.
First we show that every S-Willmore surface with flat normal bundle in $S^n$ must locate in some $S^3\subset S^n$, from which we characterize Clifford torus as the only non-equatorial homogeneous minimal surface in $S^n$ with flat normal bundle, which improve a result of K. Yang. Then we derived that every Willmore two sphere with flat normal bundle in $S^n$ is conformal to a minimal surface with embedded planer ends in $\mathbb{R}^3$. We also point out that for a class of Willmore tori, they have flat normal bundle if and only if they locate in some $S^3$. In the end, we show that a Willmore surface with flat normal bundle must locate in some $S^6$.\\

{\bf Keywords:} Willmore surfaces, S-Willmore surfaces, Willmore sphere, Clifford torus, flat normal bundle.\\

%{\bf PACS2010: 02.40.Hw, 02.40.Ma, 02.30.Ik}\\

{\bf MSC2010: 53A30, 53A07, 53B30}\\

\section{Introduction}

The fundamental work of Bryant on Willmore surfaces \cite{Bryant1984}, introduces several directions for the study of Willmore surfaces. The existence of dual surfaces reveals the transforming properties of Willmore surfaces, generalized by peoples in many papers, for example \cite{Ejiri1988,BFLPP,Ma2006}. The harmonicity of conformal Gauss map leads to the use of integrable systems in this field \cite{BFLPP,Helein,Xia-Shen,DWa1}.  And the classification theorem of Willmore 2-spheres are generalized by many peoples in all kinds of ways \cite{BFLPP,DWa1,Ejiri1988,Ma,Ma-W1,Mus1,Mon}.

Although there has been many interesting progress on this topic now, it is still not so clear for people to have enough knowledge on Willmore surfaces, mainly due to the complexity of Willmore equations and the lack of powerful tools. Then people began to understand Willmore surfaces with special conformal properties, see \cite{BPP,Ejiri1988,Ma-WC,Le-Pe-Pin} for instance.

In this paper we want to see what will happen if a Willmore surface is of flat normal bundle. This problem is motivated by the following considerations. Firstly, a theorem in \cite{BPP} tells that isothermic Willmore surfaces in $S^n$ can be reduced to some $S^4$. So it is natural to ask what will happen if we weaken the isothermic condition to be of flat normal bundle, which is also a conformal invariant property. Secondly, in another interesting paper \cite{Le-Pe-Pin}, it was shown that Willmore tori in $S^4$ with non trivial normal bundle (hence non-flat normal bundle), are in fact given by holomorphic data.
 This also induce us to consider Willmore surfaces of flat normal bundle (the simplest one of those of trivial normal bundle). The last reason is related with the classification problem of Willmore spheres. Willmore spheres in $S^3$ and in $S^4$ (which are all S-Willmore) have been classified, and studied in details in the papers of Bryant \cite{Bryant1984}, Ejiri \cite{Ejiri1988}, Musso \cite{Mus1}, Montiel \cite{Mon} etc. Ejiri also classified S-Willmore spheres in $S^n$ for all $n\geq3$. In a recent work, Dorfmeister and the author showed that there exist non S-Willmore Willmore spheres in $S^6$, and gave a classification of  Willmore spheres in $S^6$ by use of loop group methods. However, the  geometric meaning of the classification stays unclear now. So we want to have a look at Willmore spheres in a geometric way, which leads us to the consideration of Willmore spheres of flat normal bundle.

It is not easy to show the existence or non-existence of Willmore surfaces with flat normal bundle. An example is the work of Yang \cite{Yang} on homogeneous minimal surfaces with flat normal bundle in $S^n$, where he reduced such surfaces into some real analytic variety. This did not give a sufficient description of such surfaces. For our case, we add some further conditions. First we assume that the Willmore surfaces are S-Willmore. Then we show that they must reduce to some $S^3\subset S^n$. Noticing that minimal surfaces in $S^n$ are S-Willmore, we simplify the result of \cite{Yang}. In fact we show that Clifford torus is the only homogeneous minimal surface in $S^n$ with flat normal bundle. Second, we consider Willmore two spheres with flat normal bundle. Then by a detailed computation, we show that such surfaces will also reduce to some $S^3\subset S^n$. Then, one can see non S-Willmore Willmore spheres must have non flat normal bundle. Thirdly, we focus on surfaces in $S^{2n+1}$ derived from the Hopf bundle $\pi: S^{2n+1}\rightarrow \mathbb{C}P^n$, showing that such surfaces must locate in some $S^3$ if they are assumed of flat normal bundle. In the end, to compare with \cite{BFLPP}, similar to the treatment in \cite{BFLPP}, by use of detailed discussion on the integrable equations, we show that every Willmore surface with flat normal bundle locates in some $S^6$, which is weaken than the results in \cite{BFLPP}.

This paper is organized as follows. In Section~2  we
review the main theory of Willmore surfaces quickly and give some basic description on surfaces of flat normal bundle. Then we do with S-Willmore surfaces with flat normal bundle and Willmore spheres with flat normal bunlde respectively in Section~3 and Section~4. In Section 5 we discuss the surfaces related with Hopf bundle. Then we end the paper by giving some more discussion on the integrable equations in Secion 6.

\section{Preliminary}

We denote the Minkowski space-time $\mathbb{R}^{n}_1$ as $\mathbb{R}^{n}$ equipped with a Lorenzian metric
$<x,y>=-x_{0}y_0+x_1y_1+\cdots+x_ny_n.$
Let $\mathcal{C}^{n+1}$ be the light cone of $\mathbb{R}^{n+2}_{1}$. One can see that the projective light cone
$$
Q^{n}=\{\ [x]\in\mathbb{R}P^{n+1}\ |\ x\in C^{n+1}\setminus \{0\}
\}
$$
with the induced conformal metric, is conformally equivalent to $S^{n}$. And the conformal group of
$Q^{n}$ is exactly the orthogonal group $O(n+1,1)/\{\pm1\}$ of
$\mathbb{R}^{n+2}_1$, acting on $Q^{n}$ by
$
T([x])=[Tx],\ \forall T\in O(n+1,1).
$

Let $y:M\rightarrow S^{n}$ be a conformal immersion from a Riemann surface $M$. Let $U\subset M$ be an open subset. A local
lift of $y$ is a map $Y:U\rightarrow C^{n+1}\setminus\{0\}$ such
that $\pi\circ Y=y$. Two different local lifts differ by a scaling,
thus deriving the same conformal metric on $M$. Here we call $y$ a {\em conformal} immersion, if $\langle Y_{z},Y_{z}\rangle=0$ and $\langle
Y_{z},Y_{\bar{z}}\rangle >0$ for any local lift $Y$ and any complex
coordinate $z$ on $M$. There is a natural decomposition $M\times
\mathbb{R}^{n+2}_{1}=V\oplus V^{\perp}$, where
\begin{equation}
V={\rm Span}\{Y,{\rm Re}Y_{z},{\rm Im}Y_{z},Y_{z\bar{z}}\}
\end{equation}
is a Lorentzian rank-4 subbundle independent of the choice of $Y$
and $z$, and $V^{\perp}$ is the orthogonal complement of $V$. Note that both of $V$ and $V^{\perp}$ are conformal invariant. Their
complexifications are denoted separately as $V_{\mathbb{C}}$ and
$V^{\perp}_{\mathbb{C}}$.

Fix a local coordinate $z$. There is a local lift $Y$ satisfying
$|{\rm d}Y|^2=|{\rm d}z|^2$, called the canonical lift (with respect
to $z$). Choose a frame $\{Y,Y_{z},Y_{\bar{z}},N\}$ of
$V_{\mathbb{C}}$, where $N\in\Gamma(V)$ is uniquely determined by
\begin{equation}\label{eq-N}
\langle N,Y_{z}\rangle=\langle N,Y_{\bar{z}}\rangle=\langle
N,N\rangle=0,\langle N,Y\rangle=-1.
\end{equation}

The map $Gr:M\rightarrow Gr_{3,1}(R^{n+2}_1)$ defined by
$$p\in M\mapsto V|_p$$
is denoted
 \emph{the conformal Gauss map} of $y$ (See
also \cite{Bryant1984,BPP,Ejiri1988,Ma}).

Given frames as above, and noting that $Y_{zz}$ is orthogonal to
$Y$, $Y_{z}$ and $Y_{\bar{z}}$, there exists a complex function $s$
and a section $\kappa\in \Gamma(V_{\mathbb{C}}^{\perp})$ such that
\begin{equation}
Y_{zz}=-\frac{s}{2}Y+\kappa.
\end{equation}
This defines two basic invariants $\kappa$ and $s$ depending on
coordinates $z$, \emph{the conformal Hopf differential} and
\emph{the Schwarzian} of $y$ (for more discussion, see \cite{BPP,Ma}).
The conformal
Hopf differential plays an important role in the research of Willmore
surfaces. Direct computation shows that the conformal Gauss map $Gr$ induces a conformal-invariant
metric
\[
g:=\frac{1}{4}\langle {\rm d}G,{\rm d}G\rangle=\langle
\kappa,\bar{\kappa}\rangle|dz|^{2}
\]
on M. Note that this metric degenerates at umibilic points of $y$.  Now we define the Willmore
functional and Willmore surfaces by use of this metric.

\begin{definition} \emph{The Willmore functional} of $y$ is
defined as the area of M with respect to the metric above:
\begin{equation}
W(y):=2i\int_{M}\langle \kappa,\bar{\kappa}\rangle dz\wedge
d\bar{z}.
\end{equation}
An immersed surface $y:M\rightarrow S^{n}$ is called a
\emph{Willmore surface}, if it is a critical surface of the Willmore
functional with respect to any variation of the map $y:M\rightarrow
 S^{n}$.
\end{definition}
 To compare with the traditional definition, let $x:M\rightarrow \mathbb{R}^{n}$ be the spherical projection of $y$ into $\mathbb{R}^{n}$. Let $H,\ K$ denote the mean curvature and Gauss curvature of $x$. Then one can verify easily that
$$W(y)=W(x)=\int_M(H^2-K)dM,$$
coinciding with the original definition of Willmore functional.

 Let $D$
denote the normal connection and $\psi\in
\Gamma(V_{\mathbb{C}}^{\perp})$ any section of the normal bundle.
The structure equations are given as follows:
\begin{equation}\label{eq-moving}
\left\{\begin {array}{lllll}
Y_{zz}=-\frac{s}{2}Y+\kappa,\\
Y_{z\bar{z}}=-\langle \kappa,\bar\kappa\rangle Y+\frac{1}{2}N,\\
N_{z}=-2\langle \kappa,\bar\kappa\rangle Y_{z}-sY_{\bar{z}}+2D_{\bar{z}}\kappa,\\
\psi_{z}=D_{z}\psi+2\langle \psi,D_{\bar{z}}\kappa\rangle Y-2\langle
\psi,\kappa\rangle Y_{\bar{z}},
\end {array}\right.
\end{equation}
The conformal Gauss, Codazzi and Ricci equations as integrable
conditions are:
\begin{equation}\label{eq-integ}
\left\{\begin {array}{lllll} \frac{1}{2}s_{\bar{z}}=3\langle
\kappa,D_z\bar\kappa\rangle +\langle D_z\kappa,\bar\kappa\rangle,\\
{\rm Im}(D_{\bar{z}}D_{\bar{z}}\kappa+\frac{\bar{s}}{2}\kappa)=0,\\
R^{D}_{\bar{z}z}\psi=D_{\bar{z}}D_{z}\psi-D_{z}D_{\bar{z}}\psi =
2\langle \psi,\kappa\rangle\bar{\kappa}- 2\langle
\psi,\bar{\kappa}\rangle\kappa.
\end {array}\right.
\end{equation}

It is well-known that Willmore surfaces are characterized as
follow
\cite{Bryant1984,BPP,Ejiri1988,Wang1998}.

\begin{theorem}\label{thm-willmore} For a conformal immersion $y:M\rightarrow  S^{n+2}$, the following three conditions
are equivalent:

(i) $y$ is Willmore.

(ii) The conformal Gauss map $G$ is a harmonic map into
$G_{3,1}(\mathbb{R}^{n+3}_{1})$.

(iii) The conformal Hopf differential $\kappa$ of $y$ satisfies the
Willmore condition as below, which is stronger than the conformal
Codazzi equation \eqref{eq-integ}:
\begin{equation}\label{eq-willmore}
D_{\bar{z}}D_{\bar{z}}\kappa+\frac{\bar{s}}{2}\kappa=0.
\end{equation}
\end{theorem}

There is a special type of Willmore surfaces called {\em S-Willmore surfaces} first introduced by Ejiri in \cite{Ejiri1988} as the ones with dual surfaces.
\begin{definition}
A Willmore immersion $y:M^2\rightarrow S^n$ is called an S-Willlmore surface if its conformal Hopf differential satisfies
$$D_{\bar{z}}\kappa || \kappa,$$
i.e. there exists some function $\mu$ on $M$ such that $D_{\bar{z}}\kappa+\mu\kappa=0$.\\
\end{definition}

We also need the following lemma on surfaces in $S^n$ with flat normal bundle.

\begin{lemma}\label{w-nb}Let $y:M\rightarrow S^n$ be an conformal immersion, with all conformal data as above. Suppose that $y$ has no umbilic points. Then $x$ is of flat normal bundle if and only if there exists an orthonormal basis $\{\psi_{\alpha}\},\ \alpha=3,\cdots,n$ of $V^{\perp}$, such that
\begin{equation}\kappa=k_3\psi_{3}.\end{equation}
\end{lemma}
\begin{proof}
Let $\{\hat\psi_{\alpha}\},\ \alpha=3,\cdots,n$ be an orthonormal basis of $V^{\perp}$. Assume that
$$\kappa=\sum_{\alpha}\hat{k}_{\alpha}\hat\psi_{\alpha}.$$
Since $y$ has flat normal bundle, we have that
$$0=R^{D}_{\bar{z}z}\psi=
2\langle \psi,\kappa\rangle\bar{\kappa}- 2\langle
\psi,\bar{\kappa}\rangle\kappa, \forall \psi\in \Gamma(V^{\perp}).$$
So for any $\hat\psi_{\alpha}$,
$$
\langle \hat\psi_{\alpha},\kappa\rangle\bar{\kappa}- 2\langle
\hat\psi_{\alpha},\bar{\kappa}\rangle\kappa=\sum_{\beta}(\hat{k}_{\alpha}\bar{\hat{k}}_{\beta}-\bar{\hat{k}}_{\alpha}\hat{k}_{\beta})\hat\psi_{\beta}=0,$$
forcing
$$\hat{k}_{\alpha}\bar{\hat{k}}_{\beta}-\bar{\hat{k}}_{\alpha}\hat{k}_{\beta}=0, \ \forall \beta.$$
Since $y$ is umbilic free, $\sum|{\hat{k}}_{\beta}|^2>0$, we may assume that $\hat{k}_{\alpha}=|\hat{k}_{\alpha}|e^{i\theta}\neq0$, $\theta\in\mathbb{R}$. Then from above equation, we have that
$$\hat{k}_{\beta}=|\hat{k}_{\beta}|e^{i\theta}, \ \forall \beta.$$
So
$$\kappa=e^{i\theta}\hat\kappa, \hat\kappa\in \Gamma(V^{\perp}).$$
Choose a new frame  $\{\psi_{\alpha}\},\ \alpha=3,\cdots,n$ of $V^{\perp}$, with $\psi_3=\frac{1}{|\hat\kappa|}\hat\kappa$, we finish the proof.
\end{proof}

We  recall the definition of isothermic surfaces under these moving frames, see \cite{BPP,Ma} for more details.

\begin{definition}
A conformal immersion $y:M\rightarrow S^n$ is called an isothermic surface if there exist some coordinate $z$ and canonical lift $Y$ such that the corresponding conformal Hopf differential  is real, i.e.,
$\bar\kappa=\kappa.$
\end{definition}

The condition that $\kappa$ is real is equivalent to the condition that $\kappa=f\tilde\kappa$, with $\tilde\kappa$ real and $f$ holomorphic function, see \cite{BFLPP,Jer,Ma,Wang2012}.

\begin{remark}As to the difference between surfaces with flat normal bundle and isothermic surfaces, let $k_3=|k_3|e^{i\theta}$. Then by use of the elementary properties of holomorphic functions, $x$ is isothermic if and only if $\theta$ is a harmonic function, i.e., $\theta_{z\bar{z}}=0.$
\end{remark}

\section{S-Willmore surfaces with flat normal bundle}

\begin{theorem}Let $y:M\rightarrow S^n$ be an S-Willlmore surface with flat normal bundle. Then there exists a 3-dimensional sphere $S^3\subset S^n$ such that $y(M)\subset S^3$.
\end{theorem}
\begin{proof}Since $y$ is Willmore,  it is totally umbilic or has an open dense subset $M_0\subset M$ without umbilic points (This has been shown in \cite{Bryant1984} when $n=3$. A complete version for $S^n$ one may consult Proposition 5 in \cite{Ri1}). From Proposition 3.1, on $M_0$, we assume that $\kappa=k_3\psi_3$. The S-Willmore condition reads
$$D_{\bar{z}}\kappa\parallel\kappa,$$
forcing
$$D_{\bar{z}}\psi_{3}=0, \hbox{ on } M_0,$$
which, together with the structure equation, shows that $y(M_0)$ is in some three dimensional subspace $S^3$ of $S^n$. Since $M_0\subset M$ is open and dense, and Willmore surfaces are real-analytical, we have that $y(M)\subset S^3$.
\end{proof}
\begin{remark} There exist Willmore tori in $S^3$ (hence with flat normal bundle), which is not isothermic, see \cite{Pinkall1985}, \cite{FP}.

There are also many examples of Willmore tori full in $S^{2n+1}$, $n>1$, see \cite{BGS}, \cite{CF} and \cite{Li-V} for instance. However, all of these examples have non-trivial (hence non-flat) normal bundle. See Section 5 for more discussion.
\end{remark}
Since minimal surfaces in $S^n$, $H^n$ or $\mathbb{R}^n$ are S-Willmore, we have that
\begin{corollary}Minimal surfaces in  $S^n$$(H^n$ or $\mathbb{R}^n)$  with flat normal bundle must be contained in some $3$ dimensional subspace $S^3$$(H^3$ or $\mathbb{R}^3)$ of  $S^n$$(H^n$ or $\mathbb{R}^n)$.
\end{corollary}

In \cite{Yang}, Yang showed that the totality of non-equatorial homogeneous minimal surfaces in $S^{2+p}$ with
the flat normal connections is parametrized by some real algebraic variety $V \subset \mathbb{R}^{p^2}$. Here we can show that this real algebraic variety contains in fact only one point. To be concrete, we have that

\begin{theorem}  All non-equatorial homogeneous minimal surfaces in $S^{n}$ with
the flat normal connections are isometric to the Clifford torus in some $S^3\subset S^n$.
\end{theorem}

\begin{proof} The proof  comes from the fact that Clifford torus is the only non-totally geodesic minimal surface with constant Maurer-Cartan form under some parametrization.
\end{proof}
\section{Willmore spheres with flat normal bundle}

First we recall an important global holomorphic differential form defined for a Willmore surface, which is essential in the classification of Willmore spheres. For the details, we refer to \cite{Bryant1984,Ejiri1988,Ma,Ma-W1,Mon,Mus1}. Here we use the definition from Ma \cite{Ma}.
\begin{lemma}\cite{Ma} Let $y:M\rightarrow S^n$ be a Willlmore surface. Then
\begin{equation}\Omega dz^6:=(\langle D_{\bar{z}}\kappa,\kappa\rangle^2-\langle D_{\bar{z}}\kappa, D_{\bar{z}}\kappa\rangle\langle\kappa,\kappa\rangle)dz^6
\end{equation}
is a global defined holomorphic $6-$form on $M$.
\end{lemma}

\begin{theorem}Let $y:S^2\rightarrow S^n$ be a Willlmore sphere with flat normal bundle. Then $y$ is conformal to a complete, genus $0$ minimal surface in $\mathbb{R}^3$, with embedded flat ends.
\end{theorem}

\begin{proof} Since $M=S^2$, by the Riemann-Roch theorem, all holomorphic forms vanish on $S^2$, we have that
$$(\langle D_{\bar{z}}\kappa,\kappa\rangle^2-\langle D_{\bar{z}}\kappa, D_{\bar{z}}\kappa\rangle\langle\kappa,\kappa\rangle)=0.$$
Let $\{ \psi_{\alpha}\},\ \alpha=3,\cdots,n$ be an orthonormal basis of $V^{\perp}$ such that
$$\kappa= {k}_{3} \psi_{3}.$$
Then
$$D_{\bar{z}}\kappa=k_{3\bar{z}}\psi_3+k_3D_{\bar{z}}\psi_3.
$$
So
$$|k_3|^2|k_{3\bar{z}}|^2-(|k_{3\bar{z}}|^2+|k_3|^2  \langle D_{\bar{z}}\psi_3, D_{\bar{z}}\psi_3\rangle)|k_3|^2=0.$$
forcing
$$|k_3|^4 \langle D_{\bar{z}}\psi_3, D_{\bar{z}}\psi_3\rangle=0.$$
Since $y$ is totally umbilic or $k_3\neq0$ on an open dense subset $M_0\subset S^2$, we have that $y$ is totally umbilic or
$$\langle D_{\bar{z}}\psi_3, D_{\bar{z}}\psi_3\rangle=0.$$
If $y$ is totally umbilic, it is conformal to a 2-plane.
Now we suppose that \begin{equation}\label{eq-dbarz}\langle D_{\bar{z}}\psi_3, D_{\bar{z}}\psi_3\rangle=0.\end{equation}
For the same reason, just as Proposition 5 in \cite{Ri1},
$D_{\bar{z}}\psi_3\equiv0$ or $D_{\bar{z}}\psi_3\neq0$ on an open dense subset $M_1$ of $S^2$. For the first case, we have that $y$ is contained in some
$S^3\subset S^n$. For the second case, we assume that
$$D_{\bar{z}}\psi_3=b_{34}\bar{ E}_4$$
with
$$\ b_{34}=\bar{b}_{34},\  E_4\in \Gamma(V^{\perp}\otimes \mathbb{C}),\ \langle \psi_3, E_4\rangle=\langle \psi_3, \bar{E}_4\rangle=\langle E_4, E_4\rangle=0,\ \langle E_4, \bar{E}_4\rangle=2.$$
Suppose that $$D_{z} E_4= b_{44}E_4-2b_{34}\psi_3+\cdots.$$
So $$D_{z} \bar{E}_4= -b_{44}\bar{E}_4-2b_{34}\psi_3+\cdots.$$
Then Willmore equation reads
\begin{equation*}\left\{\begin{split}
                   &k_{3\bar{z}\bar{z}}+\frac{\bar{s}}{2}k_3=0, \\
                   &2k_{3\bar{z}}b_{34}+k_{3}b_{34\bar{z}}+k_{3}b_{34}\bar{b}_{44}=0.\\
                 \end{split}\right.
\end{equation*}
And the flatness of normal bundle yields
$$0=D_{\bar{z}}D_{z}\psi_3-D_{z}D_{\bar{z}}\psi_3\Longrightarrow b_{34\bar{z}}=b_{34}\bar{b}_{44}.$$
So we have that
$$2k_{3\bar{z}}b_{34}+2k_{3}b_{34\bar{z}}=2(k_{3}b_{34})_{\bar{z}}=0.$$
Since $b_{34}\neq0$ on $M_1$ and $b_{34}=\bar{b}_{34}$, one may assume that $b_{34}>0$. Now let $w=z(w)$ satisfy $$\left(\frac{\partial z}{\partial w}\right)^2=\frac{1}{b_{34}k_3}.$$
Let $$\hat{Y}=\frac{1}{|z_{w}|}Y.$$
Then $\hat{Y}$ is the canonical lift with respect to $w$ and
$$\hat\kappa=\hat{Y}_{ww}+\frac{\hat{s}}{2}\hat{Y}=\frac{1}{b_{34}|z_{w}|}\psi_3,\ \Rightarrow\ \hat{\kappa}=\bar{\hat{\kappa}}.$$
So $y$ is isothermic. By Theorem 3.4 in \cite{BPP}, we reduce $y$ into some $S^4$ of $S^n$. For \eqref{eq-dbarz} to be satisfied in $S^4$, one must have that $D_{\bar{z}}\psi_3\equiv0$, hence $y$ is contained in some $S^3$. By the classical theorem of Bryant in \cite{Bryant1984}, $y$ is conformal to some complete genus 0 minimal surface in $\mathbb{R}^3$ with embedded flat ends.
\end{proof}

\section{Willmore surfaces in $S^{2n+1}$ derived from the Hopf bundle}

 Since the work of Pinkall in \cite{Pinkall1985}, many examples of Willmore tori are obtained from the Hopf bundle $\pi:S^{2n+1}\rightarrow \mathbb{C}P^n$, see \cite{BGS}, \cite{CF}. In this section, we will show that all such surfaces must locate in some $S^3$ if they are of flat normal bundle.

\begin{proposition}Let $\gamma=\gamma(t): \mathbb{R}\rightarrow \mathbb{C}P^n$ be a regular curve in $\mathbb{C}P^n$. Let $x=\pi^{-1}\gamma:\mathbb{R}\times S^1\rightarrow S^{2n+1}$ be the $S^1-$invariant surface in $S^{2n+1}$ derived by $\gamma$ via the Hopf bundle. If $x$ is of flat normal bundle, by some conformal transform, $x$ can be located into some $S^3\subset S^{2n+1}$ and $\gamma$ is in some $\mathbb{C}P^1\subset\mathbb{C}P^n$. So if we
further assume that $x$ is Willmore, $x$ must be one of the examples given in \cite{Pinkall1985}.
\end{proposition}

\begin{proof} Let $\mathbb{C}^{n+1}$ be the $n-$dimensional complex vector space, with inner product
$$\langle Z, W \rangle=\sum_{j=1}^{n+1}z_j\bar{w}_j,\ Z=(z_1,\cdots,z_{n+1})\in \mathbb{C}^{n+1},
\ W=(w_1,\cdots,w_{n+1})\in \mathbb{C}^{n+1}. $$
Then $S^{2n+1}=\{Z\in \mathbb{C}^{n+1}|\langle Z,Z\rangle=1\}$ is the standard sphere of constant curvature one and  $\pi:S^{2n+1}\rightarrow \mathbb{C}P^n$ is the well-known Hopf bundle with fiber $S^1$. It is well-known that every fiber $S^1$ is totally geodesic in $S^{2n+1}$.

Note that although $\mathbb{C}^{n+1}$ is complex space, here it should be looked as a $2n+2$ dimensional real space. So we use $\sqrt{-1}$  to be the imaginary unit used in $\mathbb{C}^{n+1}$ and use $i$ to denote the complex structure of $M$ as above.

It is easy to check that for any $Z\in\mathbb{C}^{n+1}$, $\sqrt{-1}Z\in \mathbb{C}^{n+1}$ is perpendicular to $Z$.

Since $\gamma$ is a regular curve in $\mathbb{C}P^n$, there exists a lift $\tilde\gamma$ of $\gamma$, $\tilde\gamma:\mathbb{R}\rightarrow S^{2n+1}\subset\mathbb{C}^n$ such that $$\tilde\gamma_t=\frac{d}{dt}\tilde\gamma\perp \sqrt{-1}\tilde\gamma.$$
We also assume that $t$ is an arc parameter, that is $\langle\tilde\gamma_t, \tilde\gamma_t \rangle=1$. Let $\xi=\tilde\gamma_t$. Then we have that
$$\xi_t=k_1\sqrt{-1}\xi+k_2\eta-\tilde\gamma,$$
with $\eta$ some unit section perpendicular to $\tilde{\gamma},\  \sqrt{-1}\tilde\gamma,\ \xi,\ \sqrt{-1}\xi$.

 Now $x$ can be written as $x=e^{\sqrt{-1}\theta}\tilde\gamma$ and $z=t+i\theta$ is a complex coordinate of $x$. So
 \begin{equation*}\left\{\begin{split}
&x_t=e^{\sqrt{-1}\theta}\xi,\\
&x_\theta=e^{\sqrt{-1}\theta}\sqrt{-1}\tilde\gamma,\\
&x_{\theta\theta}=-e^{\sqrt{-1}\theta}\tilde\gamma,\\
&x_{tt}=k_1e^{\sqrt{-1}\theta}\sqrt{-1}\xi+k_2e^{\sqrt{-1}\theta}\eta-e^{\sqrt{-1}\theta}\tilde\gamma,\\
&x_{t\theta}=e^{\sqrt{-1}\theta}\sqrt{-1}\xi.\\
\end{split}\right.
\end{equation*}
Then the vector-valued Hopf differential $\Omega dz^2:=(x_{zz}\ \hbox{mod } \{x_t,x_{\theta}\})dz^2$ of $x$ is
$$\frac{1}{4}\left(k_1e^{\sqrt{-1}\theta}\sqrt{-1}\xi+k_2e^{\sqrt{-1}\theta}\eta-2ie^{\sqrt{-1}\theta}\sqrt{-1}\xi\right).$$
From the relation between $\Omega$ and the conformal Hopf differential $\kappa$ of $x$, see \cite{Wang2012}, together with Lemma \ref{w-nb}, we see that $x$ is of flat normal bundle if and only if $k_1e^{\sqrt{-1}\theta}\sqrt{-1}\xi+k_2e^{\sqrt{-1}\theta}\eta$ and $\sqrt{-1}e^{\sqrt{-1}\theta}\xi$ are linear dependent, i.e. $k_2\equiv0$. That is, $\tilde\gamma$ is located in some $\mathbb{C}^2\subset \mathbb{C}^{n+1}$. Then the proposition follows directly.
\end{proof}

\begin{remark}
From above, by directly computation, we have that
$$W(x)=\int_M(\frac{k_1^2+k_2^2}{4}+1)dtd\theta.$$
Note there has been enough descriptions on $\gamma$ in \cite{BGS} and \cite{CF}. Hence we omit further computations here and just focus on the ones of flat normal bundle.

Similar discussion shows that the examples in \cite{Li-V} can not have flat normal bundle. Here we omit it.

As to Willmore tori of flat normal bundle, there are a kind of interesting examples in \cite{Ma-W1}, which are isothermic Willmore tori in 4-dimensional Lorenzian space $Q^4_1$. It stays an interesting question that whether are there similar isothermic Willmore tori full in $S^4$.
\end{remark}

\section{Further remarks}

 It remains an open problem that whether there are Willmore surfaces of flat normal bundle and full in $S^n, n>4$. Note that for Willmore surface with flat normal bundle, the conformal Hopf differential may not be real. Then different from the isothermic ones, one can not reduce the Willmore equation to $(n-2)$ equations on a real vector bundle $M\times \mathbb{R}^{n-2}$. So the method used in \cite{BFLPP} fails partially in our case. However, one still can reduce the dimension $n$ to $6$ following the treatment in \cite{BFLPP}. To be concrete, we have that
\begin{proposition}Let $x:M\rightarrow S^n$, $n\geq6$, be a Willmore surface with flat normal bundle. Then $x$ locates in some $S^6\subset S^n$.
\end{proposition}
\begin{proof} Let $\tilde\psi$ be a parallel orthonormal basis of the normal bundle. From Lemma \ref{w-nb}, we have that
$$\kappa=e^{i\theta}\sum_{\alpha=3}^{n}\tilde{k}_{\alpha}\tilde\psi_{\alpha}=e^{i\theta}\tilde{\kappa},\ \tilde{k}_{\alpha}=\bar{\tilde{k}}_{\alpha}.$$
Together with the Willmore equation, we derive that
\begin{equation*}\left\{\begin{split}
&D_{\bar{z}}\kappa\in\hbox{Span}_{\mathbb{C}}\{\tilde{\kappa},\ D_{\bar{z}}\tilde{\kappa}\},\\
 &D_{z}\kappa\in\hbox{Span}_{\mathbb{C}}\{\tilde{\kappa},\ D_{z}\tilde{\kappa}\},\\ &D_{\bar{z}}\tilde{\kappa}\in\hbox{Span}_{\mathbb{C}}\{Re(D_{z}\tilde{\kappa}),\ Im( D_{z}\tilde{\kappa})\},\\
&D_{\bar{z}}D_{\bar{z}}\kappa\in\hbox{Span}_{\mathbb{C}}\{\tilde{\kappa}\}.
\end{split}\right.
\end{equation*}
Repeating using the Willmore equation and the flatness of normal bundle, we have that all the derivatives of $\kappa$ locates in the subspace
$$\hbox{Span}_{\mathbb{C}}\{\tilde\kappa,Re(D_{z}\tilde{\kappa}),\ Im( D_{z}\tilde{\kappa}), D_{\bar{z}}D_z\tilde\kappa\}.$$
Since
$\{\tilde\kappa,Re(D_{z}\tilde{\kappa}),\ Im( D_{z}\tilde{\kappa}), D_{\bar{z}}D_z\tilde\kappa\}$ is an real subspace of at most four dimension, $x$ is located in some $S^6\subset S^n$.
\end{proof}

\begin{remark} It remains an open problem that whether there are Willmore surfaces with flat normal bundle full in $S^6$. In this case, the integrable conditions are as follows:
\begin{equation*}\left\{\begin{split}
&\tilde{k}_{\alpha\bar{z}\bar{z}}+2i\theta_{\bar{z}}\tilde{k}_{\alpha\bar{z}}+(i\theta_{\bar{z}\bar{z}}-\theta_{\bar{z}}^{2} +\frac{\bar{s}}{2})\tilde{k}_{\alpha}=0,\ \alpha=3,4,5,6,\\
&\frac{1}{2}s_{\bar{z}}=2\left(\sum_{\alpha=3}^{6}\tilde{k}_{\alpha}^2\right)_{z}-2i\theta_{z}\sum_{\alpha=3}^{6}\tilde{k}_{\alpha}^2.\\
\end{split}\right.
\end{equation*}
Note that in this case, $
\tilde{k}_{\alpha\bar{z}\bar{z}}$ depends on $
\tilde{k}_{\alpha\bar{z}}$ and $\tilde{k}_{\alpha}$, so the further discussion in \cite{BFLPP} is invalid here.
It remains unclear that whether the above integrable equation has a nontrivial solution or not. Note that if $\tilde{k}_{3}$ and  $\tilde{k}_{4}$
are given, one can solve $\theta_{z}$ and $s$ from the Willmore equation. Then also from the Willmore equation,  $\tilde{k}_{5}$ and  $\tilde{k}_{6}$ are solved. So in general, such surfaces are determined by a pair of real functions  $\{\tilde{k}_{3}, \tilde{k}_{4}\}$ and two initial conditions $\tilde{k}_{5}(z_0)=c_5$, $\tilde{k}_{6}(z_0)=c_6$ for $\tilde{k}_{5}$ and $ \tilde{k}_{6}$. To fit with the conformal Gauss equation,  $\{\tilde{k}_{3}, \tilde{k}_{4}\}$, $c_5$ and $c_6$ should also satisfy some equation. So it remains an interesting question that whether are there Willmore surfaces full in $S^6$, of flat normal bundle.\\
\end{remark}

{\bf Acknowledgments}\\

Part of this work was done when the author visited the Department of Mathematics of Tuebingen University. He would like to express his gratitude for the hospitality and financial support. The author is grateful to Professor F. Pedit and Professor J. Dorfmeister for valuable discussions on Willmore surfaces and pointing to the author the paper of K. Yang.
The author is also thankful to Professor Xiang Ma for valuable discussions and suggestions.
\def\refname{Reference}

\vspace{5mm}  \noindent Peng Wang\\ {\em Department of
Mathematics, Tongji University, \\Siping Road 1239,  Shanghai, 200092, People's Republic
of China.}\\ E-mail: {\sf netwangpeng@tongji.edu.cn}
\end{document}